# A REVISED ANALYSIS OF THE EXIT PROBLEM AT WEAK NOISE, AND A SIMPLER COMPUTATION OF THE QUASIPOTENTIAL WITH TWO VARIABLES

Dietrich Ryter

RyterDM@gawnet.ch

Midartweg 3   CH-4500 Solothurn  Switzerland

Phone  +4132 621 13 07

A new approach for the weak noise analysis of exit problems removes an intrinsic contradiction of an existing method. It applies for both the mean time and the location of the exits; novel outcomes mainly concern the exits from an entire domain of attraction.

Moreover, the involved quasipotential is obtained without use of a Hamiltonian system in the case of two variables.

PACS :  05.10.Gg

*Key words:*  Weak noise; exit problem; quasipotential;

# I. Introduction

Noise-induced exits from a domain (in presence of a dynamical system) are described by the Feynman-Kac theory [1] . Its application in practical problems is however so intricate that alternative methods were presented, especially for weak noise and for domains with an attractor A of the dynamical system. With unit diffusion and with a gradient drift field this was accomplished in [2] . A more general theory was put forward in [3,4] with ideas of [5], but the application to the Kramers problem [4,6] reveals some deficiencies:

i)   dependence of the rate on an inaccessible part of the potential

ii)  inconsistent approaches for the rate and for the locations [7] of the exits

iii) different methods according to the friction range.

The present paper treats the non-gradient drift fields (and a general diffusion) in a novel way, and eliminates these shortcomings. It is based on the Green's function of the Fokker-Planck operator, with sources near A . As it turns out, it differs from the existing one by the need of evaluating a backward equation in terms of an "associated" drift, instead of the actual one. This drift is obtained by reversing the conservative contribution (determined by the quasipotential), which leaves the attractor unchanged. For the final result that substitution is rather insignificant, when the part of the boundary with the most likely exits is attracted to A , and when the quasipotential is regular there. However, for exits from a whole basin of attraction, the boundary integral for the rate now accounts for the asymmetric arrivals of the (actual) trajectories near the most likely exit point. This asymmetry is implied by the twist of the associated drift against the separatrix , and it removes the flaws mentioned above.

The associated drift naturally arises when the quasipotential is evaluated by the Hamiltonian method for solving the "eikonal" or Freidlin equation, since it is given by one group of the canonical equations. A further new finding states that the Hamiltonian system is not really required in the case of two variables: the quasipotential can then be computed without an extension of the variable space, and this more direct method immediately yields the abovementioned twist responsible for the asymmetry of the exits.

It will further be argued that the existing method is rather suited for *transition* rates in bistable models, which are not always simply related with the exit rates.

The paper is organized as follows: The exposition of the background and the Green's function approach for the exit problem will be presented for $n$ variables. The remainder with $n = 2$ includes the new construction of the quasipotential (generating function $\chi$, numerical aspects, expansions at equilibrium points), followed by the asymptotic solution of the new backward equation (in standard variables, according to the rank of the diffusion matrix). Finally, the Kramers problem will be revisited, in order to exhibit the theoretical progress of the novel method, as well as to put forward some new results concerning the exit rate for any friction and the location of the exit points.

## II. Background

The autonomous dynamical system with variables $x^i$, $i = 1,...,n$

$$\dot{x}^i = a^i(\vec{x}) \tag{2.1}$$

with smooth functions $a^i$ is supposed to have at least one compact attractor A.
The negative divergence (contraction) of the drift field $\vec{a}$ will be denoted by

$$\rho(\vec{x}) := -\nabla \cdot \vec{a} \ . \tag{2.2}$$

White Gaussian noise is further supposed to interfere in such a way that the pdf $w(\vec{x},t)$ obeys a Fokker-Planck equation [8,9] with the symmetric diffusion matrix $2\varepsilon \underline{D}(\vec{x})$, with smooth elements and with $\det \underline{D} := \delta \geq 0$ :

$$w_t = \nabla \cdot (-\vec{a}w + \varepsilon \underline{D}\nabla w) = \rho w - \vec{a} \cdot \nabla w + \varepsilon \nabla \cdot (\underline{D}\nabla w) \ , \tag{2.3}$$

$\varepsilon > 0$ being the noise strength. The noise-induced drift $\varepsilon D^{ij}{}_j w$, see [8,10], is accounted for in (2.3). (Partial derivatives are marked by a mere subscript, and the summation convention is understood; matrix symbols are underlined).

A steady state solution ($w_t = 0$) is supposed to assume the form

$$w_0(\vec{x}) = N \exp[-\phi(\vec{x})/\varepsilon] \tag{2.4}$$

asymptotically for small $\varepsilon$ ($\phi$ is the "quasipotential" or "eikonal"[11-13]).

Insertion of (2.4) into (2.3) results in

$$(\vec{a} + \underline{D}\nabla\phi)\cdot\nabla\phi + \varepsilon[\rho - \nabla\cdot(\underline{D}\nabla\phi)] = 0 . \tag{2.5}$$

Clearly the "eikonal" or Freidlin equation for $\phi$

$$(\vec{a} + \underline{D}\nabla\phi)\cdot\nabla\phi = 0 \tag{2.6}$$

yields the weak noise asymptotics, and (2.4) even holds for *each* $\varepsilon > 0$, when $\phi$ also fulfils the remaining

$$\rho = \nabla\cdot(\underline{D}\nabla\phi) . \tag{2.7}$$

The auxiliary drift field

$$\vec{a} + \underline{D}\nabla\phi := \vec{a}_c \tag{2.8}$$

is sometimes called "conservative" (since $\vec{a}_c \perp \nabla\phi$ by (2.6)), and (2.7) amounts to

$$\nabla\cdot\vec{a}_c = 0. \tag{2.9}$$

With the "diffusive drift"

$$\vec{a}_d := -\underline{D}\nabla\phi \tag{2.10}$$

$\vec{a}$ is thus decomposed as $\vec{a} = \vec{a}_c + \vec{a}_d$.

The usual way to solve (2.6) for $\phi$ is to consider the Hamiltonian $H = p_i(a^i + D^{ij}p_j)$ with the momenta $p_i := \partial\phi/\partial x^i$, and to integrate

$$\dot{x}^i = \partial H/\partial p_i = a^i + 2D^{ij}p_j \tag{2.11}$$

$$\dot{p}_i = -\partial H / \partial x^i = -p_k (\partial a^k / \partial x^i + p_j \partial D^{jk} / \partial x^i) \ . \tag{2.12}$$

*The Associated System*

As (2.6) shows, a reversal of $\vec{a}_c$ ($\widetilde{\vec{a}}_c = -\vec{a}_c$) does not affect $\phi$, and (2.9) is also unchanged. The corresponding "associated" drift is

$$\widetilde{\vec{a}} := \vec{a}_d - \vec{a}_c = -(\vec{a} + 2\underline{D}\nabla\phi) \ . \tag{2.13}$$

It is easily seen from (2.6) that

$$\dot{\phi} := \vec{a} \cdot \nabla\phi = -\nabla\phi \cdot \underline{D} \nabla\phi \leq 0 \quad \text{and} \quad \widetilde{\vec{a}} \cdot \nabla\phi = \vec{a} \cdot \nabla\phi \ . \tag{2.14}$$

Therefore $\phi$ is a Lyapunov function of both $\dot{\vec{x}} = \vec{a}$ and $\dot{\vec{x}} = \widetilde{\vec{a}}$.

Mind that (2.11) amounts to

$$\dot{x}^i = -\widetilde{\vec{a}} \ . \tag{2.15}$$

The associated drift thus follows the projection of the Hamiltonian curves into the $\vec{x}$-space.

In [14] the associated system was introduced *without* the eikonal form of $w_0$, and for a general Fokker-Planck operator $L$. The main difference is that there $\vec{a}_d := w_0^{-1} \varepsilon \underline{D} \nabla w_0$ (wherever $w_0 > 0$), and $\vec{a}_c := \vec{a} - \vec{a}_d$ (and $\widetilde{\vec{a}} := \vec{a}_d - \vec{a}_c$). For the adjoint or backward operator $L^+$ it was shown there that the operator identity

$$L w_0 = w_0 \widetilde{L}^+ \tag{2.16}$$

holds in general. For systems with detailed balance (2.16) was mentioned in [9], and in the Kramers problem it entails the "outer normal expansion" of [7].

### III. The Exit Problem

It is now supposed that a domain $\Omega$ of the $\vec{x}$-space contains one attractor A, and that the boundary $\partial\Omega$ is smooth and absorbing. The problem is when and where particles driven by weak noise arrive on $\partial\Omega$ after a start near A. A natural tool for analyzing this problem is the Green's function $G(\vec{x}|\vec{x}_0)$ of the Fokker-Planck operator $L$

$$L\,G(\vec{x}|\vec{x}_0) = -\delta(\vec{x}-\vec{x}_0) \quad \text{with} \quad G=0 \quad \text{for} \quad \vec{x}\in\partial\Omega \tag{3.1}$$

and with a starting point $\vec{x}_0$ near A (here $\delta$ denotes the n-dimensional deltafunction). Since $w_0$ is concentrated near A, the relevant information is already contained in

$$\overline{G}(\vec{x}) := {}_\Omega\!\int G(\vec{x}|\vec{x}_0)\,w_0(\vec{x}_0)\,d\vec{x}_0 \tag{3.2}$$

with $w_0$ normalized on $\Omega$. From (3.1) it follows that

$$L\,\overline{G} = -w_0 \ .$$

The righthand side of (3.1) describes the insertion of one particle per unit time at $\vec{x}_0$. The mean lifetime of trajectories starting near A is therefore

$$T = {}_\Omega\!\int \overline{G}(\vec{x})\,d\vec{x} \ , \tag{3.3}$$

while the density of the exit points is given by the modulus of $\underline{D}\nabla\overline{G}$ on $\partial\Omega$.

It is now useful to introduce the function $Q(\vec{x})$ by $\overline{G} = w_0 Q$ (where $w_0 > 0$). It obeys $L(w_0 Q) = -w_0$, and by (2.16) it follows that

$$\widetilde{L}^+ Q = -1 \quad \text{and} \quad Q = 0 \quad \text{on} \quad \partial\Omega \ . \tag{3.4}$$

Mind that this involves the *associated* drift. The solution of (3.4) yields both $T$ and the location of the exits.

Weak noise allows further simplifications. Note that the steady particle flux on $\partial\Omega$ is

$$1 = \varepsilon \,{}_{\partial\Omega}\!\int \underline{D}\nabla\overline{G}\cdot d\sigma = \varepsilon \,{}_{\partial\Omega}\!\int w_0\,\underline{D}\nabla Q\cdot d\sigma \qquad (d\sigma \text{ an element of } \partial\Omega)\ .$$

Dividing this by $T$, and introducing $T^{-1}Q(\vec{x}) := q(\vec{x})$ yields

$$r_e = \varepsilon\,{}_{\partial\Omega}\!\int w_0\,\underline{D}\nabla q\cdot d\sigma \tag{3.5}$$

for the exit rate $r_e = T^{-1}$. The unnormalized density of the exit points is $w_0|\underline{D}\nabla q|$.

The method is completed by the homogeneous equation for $q$

$$\widetilde{L}^+ q = -T^{-1} \approx 0 \quad \text{and} \quad q = 0 \quad \text{on } \partial\Omega , \qquad (3.6)$$

now with the inner boundary value $q \approx 1$ near A, which follows by (3.3) with $\overline{G} = w_0 Q$. Mind that $q(\vec{x})$ is the probability to reach the attractor region from $\vec{x}$ without a visit on $\partial\Omega$, and under the associated dynamics.

These findings agree with the existing approach [3,4] up to the crucial twiddle on $L^+$: that approach (which is not based on a Green's function) just uses the actual drift $\vec{a}$ for computing $Q$ or $q$. Clearly, the methods are equivalent when $\widetilde{\vec{a}} \equiv \vec{a}$, i.e. in [2], where $\vec{a}_c = \vec{0}$.

It is remarkable that in *bistable* models with two variables the existing method (with the separatrix as $\partial\Omega$) yields the rate of the *transitions* to the other attractor region. This was shown in the Chapt. 2 of [15] without use of absorption, but rather by considering the lowest nonzero eigenvalue of $L, L^+$ and the corresponding eigenfunctions. This means that the presence of the twiddle on $L^+$ distinguishes between the exit and the transition rates. Mind that a first arrival on the separatrix is followed either by a return, or by a transition to the other attractor. When these continuations are equally probable, the exits are just twice as frequent as the transitions. An example for this special case is the Kramers problem with a smooth threshold and with an appreciable friction, as will be shown below in the Chapter VIII.

## IV. Qualitative Impacts of the Modification

Here we briefly exhibit the major changes in the planar case, with arguments to be substantiated below.

*a) Exits from an entire basin of attraction*

The drift $\vec{a}$ usually has a hyperbolic point on the boundary $\partial\Omega$ of the basin $\Omega$, and $w_0$ on $\partial\Omega$ takes a (narrow) maximum there. If $\phi$ exists there and is regular, $\tilde{\vec{a}}$ is locally hyperbolic as well. Except for the situation in [2], the local eigendirections of these drifts do not coincide ($\tilde{\vec{a}}$ does not parallel $\partial\Omega$), so that $\tilde{\vec{a}}$ points *into* $\Omega$ on one side of the fixpoint, and outwards on the other side. As a consequence, $\nabla q$ changes rapidly along $\partial\Omega$, essentially from $O(\varepsilon^{-1})$ to zero. The density of the exit points is therefore skew [16], and the integrand in (3.5) is strongly asymmetric. In the existing method, however, $\nabla q$ is $O(\varepsilon^{-1/2})$ all along $\partial\Omega$.

*b) Exits at an attracted part of $\partial\Omega$*

It is supposed that $\phi$ is regular where $w_0$ is maximum on $\partial\Omega$ (usually at a point P). Since $\partial\Omega$ is smooth, it parallels $\vec{a}_c$ at P, so that the normal drift $a_\perp \neq 0$ is determined by $\vec{a}_d$, and is thus the same for both $\tilde{\vec{a}}$ and $\vec{a}$. It will be shown below that the parallel drift $\tilde{\vec{a}}_c$ ($= -\vec{a}_c$) does not contribute in the leading order when $a_\perp \neq 0$. Therefore both $\tilde{L}^+$ and $L^+$ yield essentially the same result for $\nabla q$ (with a possible modification when the diffusion matrix is singular, due to the function $g(x)$ in the Chaps. VII, VIII below).

## V. The Quasipotential With Two Variables

5.1 *A generating function for $\phi$*

The case of two variables

$$x^1 := x \quad , \quad x^2 := y \quad ; \quad \vec{a} := (a,b)$$

allows some substantial simplifications, also for the computation of $\phi$ and of $\tilde{a}, \tilde{b}$.

The orthogonality of $\vec{a}_c$ and $\nabla\phi$, see (2.6), suggests to introduce $\chi(x,y)$ by

$$\vec{a}_c = \chi \underline{N} \nabla\phi \quad \text{with} \quad \underline{N} := \begin{pmatrix} 0 & 1 \\ -1 & 0 \end{pmatrix}. \tag{5.1}$$

By $\vec{a} = \vec{a}_c + \vec{a}_d$ it follows that

$$\vec{a} = (-\underline{D} + \chi \underline{N})\nabla\phi, \tag{5.2}$$

or, in terms of $D^{11} := A$, $D^{12} = D^{21} := B$, $D^{22} := C$ ; $\delta = AC - B^2$ :

$$a = -A\phi_x - B\phi_y + \chi\phi_y \quad , \quad b = -C\phi_y - B\phi_x - \chi\phi_x \tag{5.3}$$

(which satisfies (2.6) for any $\chi$). With the reversed sign of $\chi$ this yields $\tilde{a}, \tilde{b}$:

$$\tilde{a} = a - 2\chi\phi_y \quad , \quad \tilde{b} = b + 2\chi\phi_x. \tag{5.4}$$

The inversion of (5.2), $\nabla\phi = (-\underline{D} + \chi\underline{N})^{-1}\vec{a}$, explicitly reads

$$\phi_x = [-aC + b(B - \chi)]/(\chi^2 + \delta) \quad , \quad \phi_y = [a(B + \chi) - bA]/(\chi^2 + \delta). \tag{5.5}$$

This must be a gradient. The condition $(\phi_x)_y = (\phi_y)_x$ results in

$$(a - 2\chi\phi_y)\chi_x + (b + 2\chi\phi_x)\chi_y = \rho\chi + \nabla \cdot (\underline{D}\underline{N}\vec{a}) + \delta_x\phi_y - \delta_y\phi_x. \tag{5.6}$$

With (5.5) inserted, this only involves $\vec{a}$ and $\underline{D}$. Mind that

i) (5.6) is linear in $\chi_x, \chi_y$ and can thus be solved by characteristics in the plane;

ii) by (5.4) the lefthand side equals $\tilde{a}\chi_x + \tilde{b}\chi_y = \vec{\tilde{a}} \cdot \nabla\chi$. The characteristics thus follow and yield the associated drift $\vec{\tilde{a}}$. Since $\vec{\tilde{a}}$ has the same Lyapunov function as $\vec{a}$, the characteristics end at attractors of $\vec{a}$, which means that they can be computed starting from there, with reversed "time" and with appropriate starting conditions.

The last property allows an easy computation of $\phi$ itself along the characteristics: from (2.14) and (5.5) it follows that

$$\dot{\phi} := \tilde{a}\phi_x + \tilde{b}\phi_y = -(a^2C - 2abB + b^2A)/(\chi^2 + \delta) \quad [= -\vec{a} \cdot (-\underline{D} + \chi\underline{N})^{-1}\vec{a}] . \tag{5.7}$$

The quasipotential $\phi$ is thus indeed obtained along characteristics in the plane, and (2.12) is not invoked.

*Remarks:*

1. Once $\chi$ is known, (5.5) can also be integrated along any path. An advantage of (5.7) is that it saves the storing of $\chi$.

2. The "time" in characteristic equations is merely an auxiliary variable, and its sign can be chosen so that the integration goes "upwards" (i.e. starting from an attractor) or "downwards" (starting from a hyperbolic point or a repellor).

3. When $\delta = 0$ (diffusion in one direction only), it is easier to work with $\chi^{-1} := \psi$. Then (5.5) - (5.7) become

$$\phi_x = \psi[(-aC + bB)\psi - b)] \quad , \quad \phi_y = \psi[(aB - bA)\psi + a] \tag{5.5'}$$

$$[a + 2(aB - bA)\psi]\psi_x + [b + 2(aC - bB)\psi]\psi_y = \psi[\rho + \psi \nabla \cdot (\underline{D}\underline{N}\vec{a})] . \tag{5.6'}$$

$$\dot{\phi} = -(a^2C - 2abB + b^2A)\psi^2 . \tag{5.7'}$$

4. When $\underline{D} = c(x,y)\underline{I}$ ($\underline{I}$ := unit matrix, $c > 0$), (5.5) entails

$$|\nabla\phi|^2 = (a^2 + b^2)/(c^2 + \chi^2) , \quad \text{and thus the bound} \quad |\nabla\phi| \leq |\vec{a}|/c .$$

5.2 *Expansions at an equilibrium point*

A vanishing drift $\vec{a} = \vec{0}$ implies $\nabla\phi = \vec{0}$ by (5.5), and therefore also $\tilde{\vec{a}} = \vec{0}$ by (2.13) or (5.4). It further follows by (5.6) that there

$$\chi = [Ca_y + B(a_x - b_y) - Ab_x]/\rho . \tag{5.8}$$

The local drift matrix

$$\underline{M} := \begin{pmatrix} a_x & a_y \\ b_x & b_y \end{pmatrix} \tag{5.9}$$

plays an important role (clearly $tr \underline{M} = -\rho$). A useful expression for it follows by (5.2), which implies

$$\underline{M} = (-\underline{D} + \chi \underline{N}) \underline{S} \qquad \text{where} \qquad \underline{S} := \begin{pmatrix} \phi_{xx} & \phi_{xy} \\ \phi_{xy} & \phi_{yy} \end{pmatrix}.$$

An immediate consequence is

$$\underline{S} = (-\underline{D} + \chi \underline{N})^{-1} \underline{M} . \tag{5.10}$$

(In the special case $\delta = 0 = \chi$, $\underline{S}$ can only exist when $\det \underline{M} = 0$). We mention that (5.3) entails $A\phi_{xx} + 2B\phi_{xy} + C\phi_{yy} = \rho$, which locally satisfies (2.7).

The emphasis lies on the analogue of $\underline{M}$ for the associated drift

$$\underline{\tilde{M}} = (-\underline{D} - \chi \underline{N}) \underline{S} = (\underline{D} + \chi \underline{N})(\underline{D} - \chi \underline{N})^{-1} \underline{M} . \tag{5.11}$$

By $\det(\underline{D} \pm \chi \underline{N}) = \delta + \chi^2$ it is readily seen that $\det \underline{\tilde{M}} = \det \underline{M}$. Since also $tr \underline{\tilde{M}} = tr \underline{M}$ (to be inferred from (5.3)), the eigenvalues of both $\underline{\tilde{M}}$ and $\underline{M}$ are the same. A more explicit form of (5.11) is

$$\underline{\tilde{M}} = \underline{K} \underline{M} \quad , \quad \underline{K} := (\delta + \chi^2)^{-1}[(\delta - \chi^2)\underline{I} + 2\chi \underline{DN}] \tag{5.12}$$

with the unit matrix $\underline{I}$. For $\delta = 0$ $\underline{K}$ reduces to $-\underline{I} + 2\chi^{-1}\underline{DN}$, and for $\underline{D} = \sqrt{\delta}\underline{I}$ it is orthogonal with the angle of rotation $2\arctan(\chi\delta^{-1/2})$. Note that $\underline{\tilde{M}} = \underline{M}$ for $\chi = 0$ only.

The matrix $\underline{\tilde{M}}$ is crucial for the exit problem, when it refers to a saddlepoint (i.e. to a hyperbolic equilibrium point on the separatrix of $\vec{a}$), since it determines $\vec{\tilde{a}}$ in the linear approximation. The eigenvalues are then realvalued, with opposite signs. Clearly, the eigenvectors $\vec{e}_-$ of $\underline{M}$ and $\vec{\tilde{e}}_-$ of $\underline{\tilde{M}}$ span the respective separatrices of $\vec{a}$ and $\vec{\tilde{a}}$; these coincide when $\chi = 0$, and are twisted otherwise. Furthermore, the eigenvector $\vec{\tilde{e}}_+$ of $\underline{\tilde{M}}$ determines the *only* characteristic leaving the hyperbolic point towards

A (and also the one towards the outer attractor, possibly at infinity); these provide the respective differences in $\phi$ by a single integration, and thus the leading contribution (i.e. the Arrhenius factor) of the exit rates from there.

*Remarks:*

1. A (hyperbolic) point with $\rho = 0$ is exceptional, as (5.6) leaves $\chi$ undetermined. The eikonal form is then only possible when also $\nabla \cdot (\underline{D}\,\underline{N}\,\vec{a}) = Ca_y + B(a_x - b_y) - Ab_x = 0$, but both $\underline{S}$ and $\underline{\tilde{M}}$ remain unspecified. According to [17] it is then also possible that several characteristics link the hyperbolic point with the same attractor.

2. The corresponding expansion in the Hamiltonian approach involves a $4 \times 4$ matrix with four eigenvalues and eigenvectors, see [18].

Since a numerical integration cannot be started on an equilibrium point, but only in its neighbourhood, it is essential to evaluate the local $\nabla \chi$ as well. This can be done by taking the derivatives in $x, y$ of (5.6). With $\nabla \cdot (\underline{D}\,\underline{N}\,\vec{a}) := \kappa$ this results in

$$\begin{pmatrix} \tilde{a}_x - \rho & \tilde{b}_x \\ \tilde{a}_y & \tilde{b}_y - \rho \end{pmatrix} \begin{pmatrix} \chi_x \\ \chi_y \end{pmatrix} = \begin{pmatrix} \rho_x \chi + \kappa_x + \delta_x \phi_{xy} - \delta_y \phi_{xx} \\ \rho_y \chi + \kappa_y + \delta_x \phi_{yy} - \delta_y \phi_{xy} \end{pmatrix} , \qquad (5.13)$$

where $\underline{\tilde{M}}$ (and $\underline{S}$ if $\nabla \delta \neq \vec{0}$) is to be used. The determinant of the matrix equals $\tilde{\delta} - \rho \, tr\underline{\tilde{M}} + \rho^2 = \delta + 2\rho^2$ and is only zero when both $\delta = 0$ and $\rho = 0$ (recall the problems when $\rho = 0$).

*Remark:*

Formally the eikonal equation (2.6) admits further solutions (always $\phi = const.$; in the Kramers model it is also easy to show that an extra nonzero $\underline{S}$ exists at a saddlepoint, which allows to construct different Hamiltonian curves). Despite the equivalence of (5.1) with (2.6), the new method yields a unique solution, at least in a finite surrounding of an equilibrium point. Only this one also satisfies (2.7) at the point itself.

### 5.3 *Exact solutions*

According to (2.9) a solution is valid for *each* $\varepsilon > 0$ when $\nabla \cdot \vec{a}_c = 0$. By (5.1) this amounts to

$$\chi_x \phi_y - \chi_y \phi_x = 0 , \qquad (5.14)$$

which holds for any $\chi = \chi(\phi)$, in particular for a constant $\chi$. By (5.6) a constant $\chi$ is given by $0 \equiv \rho\chi + \nabla \cdot (\underline{D}\underline{N}\vec{a}) + \delta_x \phi_y - \delta_y \phi_x$. If $\delta$ is constant (thus for any singular $\underline{D}$, while for a regular $\underline{D}$ this can be arranged by a change of the variables [19] ), $\chi = -\rho^{-1} \nabla \cdot (\underline{D}\underline{N}\vec{a})$, which must be a constant. The condition

$$\nabla \cdot (\underline{D}\underline{N}\vec{a}) \propto \rho \qquad (5.15)$$

is then sufficient for an exact solution. For $\underline{D} = c\underline{I}$ ($c$ a constant) this amounts to

$$\operatorname{curl} \vec{a} \propto \rho , \qquad (5.16)$$

which trivially holds for a gradient drift ($\operatorname{curl} \vec{a} = 0$). A prominent example for (5.15) is also the Kramers model ([4,6] and the Chap.VIII below), even with a state-dependent friction parameter $\gamma(x,v)$.

In contrast to "detailed balance" [9] neither (5.14) nor (5.15) is based on time reversal. Clearly, $\phi$ itself is obtained by (5.5).

## VI. The Function $q(x,y)$ For Regular Diffusion

The aim is to solve (3.6) in the plane and asymptotically for small $\varepsilon$. This is conveniently done in standardized variables so that $\underline{D} \equiv \underline{I}$ when the original $\underline{D}$ is regular everywhere (see [19] for the transformation). Then (3.6) reads

$$\tilde{a}q_x + \tilde{b}q_y + \varepsilon(q_{xx} + q_{yy}) = 0 \quad , \quad q = 0 \text{ on } \partial\Omega \quad \text{and} \quad q \approx 1 \text{ near A}. \tag{6.1}$$

### 6.1 *Reduction to an ordinary equation for small $\varepsilon$*

With the arclength $s$ along $\partial\Omega$ and $r$ normal to it ($r > 0$ in $\Omega$), and with the normal and parallel components of $\tilde{a}$, the equation (6.1) becomes

$$\tilde{a}_\perp q_r + \tilde{a}_\| q_s + \varepsilon(q_{rr} + q_{ss}) = 0 \quad , \quad q(0,s) = 0 \, , \, q(r,s) \approx 1 \text{ for large } r , \tag{6.2}$$

and it is assumed that $\tilde{a}_\perp = \alpha(s) + \beta(s)r + O(r^2)$. Clearly on $\partial\Omega$ $q_s = 0 = q_{ss}$ ; for $r > 0$ these terms remain negligible for small $\varepsilon$, as can be shown *a posteriori.*
The remaining equation $(\alpha + \beta r)q_r + \varepsilon q_{rr} = 0$ is ordinary, and solved by

$$q(r,s) = c(s) \int_0^r \exp[-\varepsilon^{-1}[\alpha(s)u + \beta(s)u^2/2]] du \quad \text{with} \tag{6.3}$$

$$c(s) := 1 / \int_0^\infty \exp\{-\varepsilon^{-1}[\alpha(s)r + \beta(s)r^2/2]\} dr \tag{6.3'}$$

when $\beta > 0$, so that

$$|\nabla q| = q_r(0,s) = c(s) . \tag{6.4}$$

Clearly $c = [2\beta/(\pi\varepsilon)]^{1/2}$ where $\alpha = 0$, and for $\alpha > 0$ the asymptotic expansion of the probability integral $\Phi$ yields

$$c \approx \alpha/\varepsilon \quad \text{when} \quad \varepsilon\beta << \alpha^2/2 . \tag{6.5}$$

Considering (6.3) as the $q^*$ of the Appendix, one can easily see that it fulfills (6.2) for small enough $\varepsilon > 0$ (if $\alpha > 0$). Mind that $\tilde{a}_\|$ is thus not involved in the leading order.

*Remark:*

The mixed derivative $q_{rs}$ of (6.3) does not vanish on $\partial\Omega$, but does not occur in (6.2).

A more general $\underline{D}$ would thus only be admitted when a principal axis parallels $\partial\Omega$.

For $\beta = 0$  $c = \alpha/\varepsilon$ if $\alpha > 0$, and $c = 0$ otherwise. When $\beta < 0$ ($\alpha > 0$), the integrand in (6.3) has the minimum value $\exp[-\alpha^2/(2\varepsilon|\beta|)]$ at $u = \alpha/|\beta| := \hat{r}$, and $r$ must be confined to $r \leq \hat{r}$ (the upper limit in (6.3') becomes $\hat{r}$). With $u := \varepsilon z$ (6.3') can be rewritten as

$$c^{-1} = \varepsilon \int_0^{\hat{r}/\varepsilon} \exp(-\alpha z + \varepsilon|\beta|z^2/2)\,dz \approx \varepsilon/\alpha \quad \text{for small enough } \varepsilon. \tag{6.6}$$

This holds when the minimum value $\exp[-\alpha^2/(2\varepsilon|\beta|)] \approx 0$, i.e. when

$$\varepsilon|\beta| \ll \alpha^2/2. \tag{6.7}$$

Note the analogy with (6.5).

### 6.2 *Attracted exit region*

Here it is supposed that at the point P, where $w_0$ is maximum on $\partial\Omega$, the drift $\vec{a}$ points inwards ($a_\perp > 0$). Clearly the smooth $\partial\Omega$ is orthogonal to $\nabla\phi$ ($\neq \vec{0}$) at P, and therefore it parallels the conservative drift $\vec{a}_c$. This implies that the normal components of both $\vec{a}$ and $\tilde{\vec{a}}$ coincide: $\tilde{a}_\perp = a_\perp$. It follows that

$$|\nabla q| = a_\perp/\varepsilon \tag{6.8}$$

irrespective of $\beta$ (mind that $\alpha > 0$, so that (6.7) holds for small enough $\varepsilon$).

Since $\tilde{a}_\| = -a_\|$ is not involved, the result is the same as for the existing approach.

It is further easily seen that the density of the exit points has its maximum at P.

### 6.3 *Exits near a transition point*

It is now supposed that $\Omega$ is the whole domain attracted to A by $\vec{a}$, and that P ($= (0,0)$ say) is an equilibrium point on $\partial\Omega$, where $w_0$ on $\partial\Omega$ is maximum ($\phi$ minimum).

Near P $\partial\Omega$ coincides with the eigenvector $\vec{e}_-$ of $\underline{M}$, while $\tilde{\vec{a}}(\vec{x}) = \underline{\tilde{M}}\vec{x}$ (with a regular $\phi$). According to (5.12) (and with $\underline{D} = \underline{I}$), $\underline{\tilde{M}}$ is given by

$$\underline{\tilde{M}} = \underline{O}^2 \underline{M} \ , \quad \underline{O} = (1+\chi^2)^{-1/2} \begin{pmatrix} 1 & \chi \\ -\chi & 1 \end{pmatrix} \tag{6.9}$$

where $\underline{O}$ describes a rotation by the angle $\arctan \chi$, while by (5.8)

$$\chi = (a_y - b_x)/\rho = \operatorname{curl} \vec{a} / \operatorname{div} \vec{a} \ . \tag{6.10}$$

Let $\lambda_\pm$ be the common eigenvalues of $\underline{M}$ and $\underline{\tilde{M}}$, $\vec{e}_\pm$ the normalized eigenvectors of $\underline{M}$, furthermore $s = 0$ at P ($s > 0$ where $\vec{\tilde{a}}$ is directed inwards). Then

$$\alpha(s) = s\, \vec{e}_+ \cdot \underline{\tilde{M}} \vec{e}_- = s\lambda_- \vec{e}_+ \cdot \underline{O}^2 \vec{e}_- = s(-\lambda_-) 2|\chi|/(1+\chi^2) \tag{6.11}$$

$$\beta = \vec{e}_+ \cdot \underline{\tilde{M}} \vec{e}_+ = \lambda_+ \vec{e}_+ \underline{O}^2 \vec{e}_+ = \lambda_+ (1-\chi^2)/(1+\chi^2) \ . \tag{6.12}$$

On $\partial\Omega$ $\phi$ reads $\phi(s) = \phi_P + k s^2/2$, ($k > 0$), where by (5.10)

$$k = \vec{e}_- \cdot \underline{S} \vec{e}_- = \lambda_- \vec{e}_-(-\underline{I} + \chi \underline{N})^{-1} \vec{e}_- = -\lambda_-/(1+\chi^2) \ . \tag{6.13}$$

a) $\beta > 0$, i.e. $\chi^2 < 1$

The integral in (6.3') exists for each $s$ (near P) and yields (6.4). Together with (6.13) this specifies both the exit rate and the density of the exit points. Mind that $\nabla q$ is only symmetric in $s$ (locally) when $\chi = 0$.

b) $\beta = 0$, i.e. $\chi^2 = 1$

For $s \geq 0$ $\quad |\nabla q| = \alpha(s)/\varepsilon = (-\lambda_-)s/\varepsilon \quad$ , while $|\nabla q| = 0$ for $s \leq 0$. $\tag{6.14}$

c) $\beta < 0$, i.e. $\chi^2 > 1$

The condition (6.7) is not met for small $|s|$. Yet one may modify (6.6) by $\exp(\varepsilon|\beta|z^2/2) \approx 1 + \varepsilon|\beta|z^2/2$ (mind that $|\beta| \leq 1$), which allows $\hat{r}/\varepsilon \approx \infty$.

The results is $\varepsilon(\alpha^{-1} + \varepsilon|\beta|\alpha^{-3})$, whence

$$\varepsilon\alpha^{-1} c \approx (1+\varepsilon|\beta|/\alpha^2)^{-1} = \alpha^2/(\alpha^2 + \varepsilon|\beta|) = 1 - \varepsilon|\beta|/(\alpha^2 + \varepsilon|\beta|) \ . \tag{6.15}$$

This exhibits the impact of $|\beta|$ for small $\alpha$. Indeed, with (6.11) rewritten as $\alpha = ms$, this shows that $c \propto s^3$ for small $s$, while $c \approx \alpha/\varepsilon$ for $s^2 >> |\beta| m^{-2} \varepsilon$, in agreement with (6.7). Clearly, $c = 0$ for $s < 0$. The result for the exit rate can be written as

$$r_e = w_0(P)\varepsilon\, m[k^{-1} - |\beta| \int_0^\infty (2mu + |\beta|)^{-1} \exp(-ku)\,du] \,. \tag{6.16}$$

For $\beta \leq 0$ the density of the exit points is conveniently expressed in terms of the scaled $s' := s\varepsilon^{-1/2}$. Up to normalization it reads

$$\exp[-k(s')^2/2]\, ms'\{1 - [|\beta|^{-1}(ms')^2 + 1]^{-1}\} \quad, \quad \text{and } 0 \text{ for } s' < 0 \tag{6.17}$$

### VII. The Function $q(x,y)$ For Singular Diffusion

It is now assumed that $\underline{D}$ has rank 1 everywhere (unidirectional diffusion), and that the variables have been transformed (see again [19]) to yield

$$\underline{D} = \begin{pmatrix} 0 & 0 \\ 0 & C \end{pmatrix} \quad \text{(with a constant } C > 0 \text{)} \,. \tag{7.1}$$

Then the equation for $q$ reads

$$\widetilde{a} q_x + \widetilde{b} q_y + \varepsilon C q_{yy} = 0 \quad, \quad q = 0 \text{ on } \partial\Omega \text{ and } q \approx 1 \text{ near } A\,. \tag{7.2}$$

Again the equation can be reduced to an ordinary one (for small $\varepsilon$), now by considering the vector $(1,\mu)$ tangent to $\partial\Omega$ and by observing that $q_x + \mu q_y = 0$ on $\partial\Omega$. Close enough to $\partial\Omega$ one may use this relation to eliminate $q_x$, which gives

$$(\widetilde{b} - \mu\widetilde{a}) q_y + \varepsilon C q_{yy} = 0 \quad \text{near } \partial\Omega, \quad q = 0 \text{ on } \partial\Omega\,. \tag{7.3}$$

Clearly, this elimination of $q_x$ is only exact on $\partial\Omega$, but the ensuing approximation for $q$ can be checked by the method of the Appendix (as in the Chapt. VIII below). With

$$\alpha(x) = \widetilde{b}(x, y_{\partial\Omega}) - \mu\widetilde{a}(x, y_{\partial\Omega}) \quad \text{and} \quad \beta(x) = \widetilde{b}_y(x, y_{\partial\Omega}) - \mu\widetilde{a}_y(x, y_{\partial\Omega})\,, \tag{7.4}$$

$$[\alpha + \beta(y - y_{\partial\Omega})]q_y + \varepsilon C q_{yy} = 0 \quad , \quad q = 0 \text{ for } y = y_{\partial\Omega} \tag{7.5}$$

can be solved as in Chapt. VI . Mind that $\alpha$ is the part of $\tilde{a}$ normal to $\partial\Omega$. Yet an important novelty arises from the fact that axes parallel to the $y$-axis (along which the diffusion acts and the integration is performed) need not cross the attractor region where $q \approx 1$. The actual $\max q(x, y) := g(x)$ may thus be $< 1$. This $g(x)$ is the missing inner boundary value in (7.3), and therefore it multiplies the above solution. The evaluation of $g(x)$ is a nontrivial task, which will only be accomplished here for the Kramers model.

When the maximum point of $w_0$ on $\partial\Omega$ is *attracted*, it follows that

$$q_y(x, y_{\partial\Omega}) = (\varepsilon C)^{-1} g(x) \alpha(x) \tag{7.6}$$

in analogy with Chapt. VI . The results are thus again the same as for the existing method, up to the factor $g(x)$, which at attracted points tends to 1 when $\varepsilon \to 0$.

It is worth noting that $C$ cancels in the integrand of (3.5) since $\underline{D} \propto C$.

For the *exits on a separatrix* with a hyperbolic point (where $x = 0 = y$), the matrix $\underline{\tilde{M}} = \underline{K}\underline{M}$, which determines $\tilde{a} = \underline{\tilde{M}}\vec{x}$, is again crucial. It is essentially specified by

$$\underline{K} = -\underline{I} + 2\psi \underline{D}\underline{N} \quad , \quad \psi := \chi^{-1} , \tag{7.7}$$

which follows by (5.12) for any singular $\underline{D}$. The more specific form (7.1) entails $\psi = \rho/(Ca_y)$ by (5.8), and thereby

$$\psi \underline{D}\underline{N} = (\rho/a_y)\begin{pmatrix} 0 & 0 \\ -1 & 0 \end{pmatrix} \quad . \tag{7.8}$$

Clearly the above tangent vector $\vec{t} := (1, \mu)$ parallels the eigenvector $\vec{e}_-$ of $\underline{M}$, and $\tilde{a}, \tilde{b}$ at $x\vec{t}$ are given by $x\underline{\tilde{M}}\vec{t} = x\lambda_- \underline{K}\vec{t}$, i.e.

$$\tilde{a}(x, \mu x) = x(-\lambda_-) \quad , \quad \tilde{b}(x, \mu x) = x(-\lambda_-)[2(\rho/a_y) + \mu] \ .$$

This yields

$$\alpha(x) = x(-\lambda_-) 2\rho/a_y \ . \tag{7.9}$$

(When $a_y = 0$, the separatrix is tangent to the $y$–axis). It is further easily seen that

$$\beta = \tilde{b}_y - \mu \tilde{a}_y = \tilde{M}_{22} - \mu \tilde{M}_{12} . \tag{7.10}$$

*Remark:*

When $g(0) = 0$, it seems natural to expect that $\beta$ is not be involved, so that (7.6) applies with (7.9). In the Kramers case this will indeed be shown below, by use of the explicit $g(x)$.

Some properties of $\underline{M}$ and $\underline{\tilde{M}}$ are particularly simple when $a_x = 0$ (as in the Kramers model below): The eigenvalues and eigenvectors are then given by

$$\lambda_\pm = [b_y \pm (b_y^2 - 4a_y b_x)^{1/2}]/2 \quad ; \quad \vec{e}_\pm \propto (a_y, \lambda_\pm) \quad , \quad \tilde{\vec{e}}_\pm \propto (\tilde{a}_y, \lambda_\pm) \tag{7.11}$$

while 
$$\underline{\tilde{M}} = \begin{pmatrix} 0 & -a_y \\ -b_x & -2(\rho/a_y)a_y - b_y \end{pmatrix} = \begin{pmatrix} 0 & -a_y \\ -b_x & b_y \end{pmatrix} \ . \tag{7.12}$$

Clearly, $a_y, b_x$ stem from $\vec{a}_c$, and $b_y$ from $\vec{a}_d$. Note in particular that

$$\tilde{\vec{e}}_+ \propto (-a_y, \lambda_+) . \tag{7.13}$$

## VIII. The Kramers Problem

The very well-known Kramers problem [4,6] is revisited here in some detail. The purpose is to show how the deficiencies of the existing method (as mentioned in the Introduction) are now resolved. This includes the explicit determination of $g(x)$, and thereby a correction of the exit point density. The usual linearization at the saddlepoint will be avoided, so that the results hold for each friction $\gamma$.

The model describes a massive particle moving in a potential $U(x)$ with a threshold at $x = 0$ (supposed to be smooth; $U(0) = 0$) and with a minimum at $x_A < 0$. The Langevin equations (with unit mass) are

$$\dot{x} = v \ , \quad \dot{v} = -\gamma v - U'(x) + (2\gamma\varepsilon)^{1/2}\xi \qquad (\xi \text{ is Gaussian white noise}) .$$

With $y := v$ this corresponds to the drift $a = v$, $b = -\gamma v - U'(x)$.

Integral curves $v(x)$ of that drift obey

$$v'v = -\gamma v - U'(x) , \qquad (8.1)$$

and $v'(0) = \lambda_\mp$ determine the separatrix $\partial\Omega$ (with $\lambda_-$), as well as the unstable manifold (see (7.11) for $\lambda_\mp$, with $a_v = 1$, $b_x = -U''(0) := \omega^2$, $b_v = -\gamma$).

The separatrix $v_s(x)$ is depicted in [4] and [15] ; note that $v_s(0) = 0 = v_s(x_1)$, where $x_1 < x_A$ is the turning point, and that $v_s(x) > 0$ in between. Mere inspection of such a plot shows that a particle coming from A cannot reach a positive $x$, because at $x = 0$ it must have a negative $v = \dot{x}$ driving it back, and the diffusion only acts in the $v$-direction. The far side ($x > 0$) of $v_s(x)$ is thus strictly *inaccessible*. This contrasts with the fact that in the existing method $\nabla q$ is $O(\varepsilon^{-1/2})$ on both sides of the threshold !

The energy on the separatrix $E_s(x) = U(x) + v_s^2(x)/2$ fulfills $E_s' = -\gamma v_s$ in view of (8.1), and is thus given by

$$E_s(x) = -\gamma I(x) \quad \text{with} \quad I(x) := \int_0^x v_s(z)\, dz . \qquad (8.2)$$

In the matrix (7.1) $C = \gamma$, and the steady-state pdf $w_0(x,v)$ reads

$$N\varepsilon^{-1} \exp\{-[U(x) + v^2/2]/\varepsilon\} \quad \text{with} \quad N := (2\pi)^{-1}(U_A'')^{1/2} \exp(U_A/\varepsilon) \qquad (8.3)$$

which amounts to $\phi(x,v) = U(x) + v^2/2$. The conservative drift (2.8) is thus $\vec{a}_c = (v, -U')$, and the associated drift $\tilde{\vec{a}} = \vec{a} - 2\vec{a}_c$ becomes

$$\tilde{a} = (-v, -\gamma v + U') . \tag{8.4}$$

The separatrix of $\tilde{a}$ is given by $\tilde{v}_s(x) = -v_s(x)$. Mind that the continuation of $v_s$ beyond the turning point (to $v_s < 0$) is not attracted to A by $\tilde{a}$, since it is separated from A by $\tilde{v}_s(x)$.

The equation for $q$ is

$$\tilde{a}q_x + \tilde{b}q_v + \varepsilon\gamma q_{vv} = -vq_x + [-\gamma v + U'(x)]q_v + \varepsilon\gamma q_{vv} = 0 \quad , \quad q[x, v_s(x)] = 0 . \tag{8.5}$$

According to (7.4) with $\mu = v_s{'}$

$$\alpha(x) = \tilde{b}(x, y_{\partial\Omega}) - \mu\tilde{a}(x, y_{\partial\Omega}) = -\gamma v_s + U' + v_s{'}v_s = -2\gamma v_s(x) \tag{8.6}$$

in view of (8.1). The negative sign of $\alpha$ is due to the fact that $v < v_s$ in $\Omega$.

The *function* $g(x)$ is the probability to reach the region of the attractor A from $x \leq 0$ in the associated dynamics, more precisely the maximum probability with respect to $v$. At $x = 0$ the nonabsorbed trajectories have a negative $v$ and continue to $x > 0$ by $\dot{x} = -v$ (no diffusion in the $x$- direction), so that $g(0) = 0$. For $x << 0$ $g \approx 1$.

The equation for $g(x)$ will be established with the actual drift ($\dot{x} = -v$ is very unintuitive), and the simultaneous validity with $\tilde{a}$ will be clear afterwards. It is based on the finding that near the top of a *smooth* threshold, and for small $|v|$, the motion is approximately Markovian in $x$ alone. Indeed, a noiseless trajectory spends an unbounded time there when it returns arbitrarily close to the top, even when $\gamma = 0$ (passing ones are absorbed). On the other hand, for $\varepsilon > 0$ (and near the top) $v(t)$ is an Ornstein process with variance $\varepsilon$ and correlation time $\gamma^{-1}$. For small enough $|x|$ and $\gamma > 0$ the process $v(t)$ has thus time enough to relax to the value $\lambda_+ x$ on the unstable manifold. This eliminates $v$ from the process $x(t)$, which is thus a "creeping" motion with the drift $\lambda_+ x$ (for each $\gamma$, but for small $|x|$ only). It obeys the Fokker- Planck equation

$$w_t = -(\lambda_+ x w)_x + \varepsilon |\lambda_-|^{-1} w_{xx}$$

where the diffusion parameter is specified to yield the stationary $w_0(x) = \exp[-U(x)/\varepsilon]$ $\approx \exp[\omega^2 x^2 (2\varepsilon)^{-1}]$ (mind that $\lambda_+ \lambda_- = -\omega^2$). [Only for large $\gamma$ this becomes the Smoluchowski equation [4], since then $|\lambda_-| \to \gamma$].

The function $g(x)$ is determined by the corresponding backward equation

$$\lambda_+ x g' + \varepsilon |\lambda_-|^{-1} g'' = 0 \text{ , or}$$

$$\omega^2 x g' + \varepsilon g'' = 0 \text{ , with } g(0) = 0 \text{ and } g(x_A) \approx 1 \text{ .}$$

The result is

$$g'(x) = -2\omega(2\pi\varepsilon)^{-1/2} \exp[-\omega^2 x^2 (2\varepsilon)^{-1}] \quad \text{and} \quad g(x) = \Phi[(2\varepsilon)^{-1/2} \omega |x|] \text{ .} \tag{8.7}$$

Note that this does not depend on $\gamma$. The simultaneous validity for the associated system is due to the fact that the only parameters are $\omega$ and $\varepsilon$ (and it can be verified more formally).

The pertinent approximation for $q$ is now

$$q(x,\eta) = g(x)\{1 - \exp[-\alpha(x)(\varepsilon\gamma)^{-1}\eta]\} \quad , \quad \eta := v - v_s(x) \leq 0$$

$$= g(x)\{1 - \exp[2v_s(x)\varepsilon^{-1}\eta]\} \text{ .} \tag{8.8}$$

It depends on $\gamma$ via $v_s$ only, and it will be confirmed in the Appendix. Accordingly

$$q_y(x) = g(x)(\varepsilon\gamma)^{-1}\alpha(x) = -2g(x)v_s(x)/\varepsilon \quad \text{on the separatrix.} \tag{8.9}$$

*Result for the exit points*

It is natural to consider the density with respect to $x$ (the arclength on the separatrix involves the sum of terms with different units). Up to normalization it is given by $|\alpha(x)| g(x) \exp[-E_s(x)/\varepsilon]$, thus by

$$v_s(x)\Phi[(2\varepsilon)^{-1/2}\omega|x|]\exp[\gamma I(x)/\varepsilon] \quad , \quad x_1 \leq x \leq 0. \tag{8.10}$$

Recall that $I(x) < 0$ for $x < 0$. The normalizing prefactor is

$$\gamma\varepsilon^{-1}\{\omega|\lambda_-|^{-1} - \exp[\gamma\varepsilon^{-1}I(x_1)]\}^{-1} \quad , \tag{8.11}$$

see the analysis for the rate below. The result of [7] would be recovered with $g(x) \equiv 1$. The actual $g(x)$ shifts the exit points away from the saddle, roughly by $(2\varepsilon)^{1/2}/\omega$.

In the limit $\gamma \to 0$ (8.10) becomes $[-U(x)]^{1/2} g(x)$.

*The exit rate*

In view of (3.5) the rate $r_e$ is now given by

$$r_e = \varepsilon\gamma\int_0^{x_1} g(x)\alpha(x)(\varepsilon\gamma)^{-1} w_0(0,0)\exp[-E_s(x)/\varepsilon]dx$$

$$= -\varepsilon^{-1}N\exp(U_A/\varepsilon)\int_0^{x_1} g(x)2\gamma v_s(x)\exp[\gamma I(x)/\varepsilon]dx$$

(recall that $x_1 < 0$). By (8.2) $v_s = I'$, so that $[\exp(\gamma I/\varepsilon)]' = \exp(\gamma I/\varepsilon)\gamma v_s/\varepsilon$.

Therefore

$$r_e = -2N\exp(U_A/\varepsilon)\int_0^{x_1} g(x)\{\exp[\gamma I(x)/\varepsilon]\}'\, dx \; .$$

Integration by parts (with $g(0) = 0$ and $g(x_1) = 1$) yields

$$r_e = 2N\exp(U_A/\varepsilon)(R - Z) \quad , \quad \text{where} \tag{8.12}$$

$$R := \int_0^{x_1} g'(x)\exp[\gamma I(x)/\varepsilon]dx \quad , \quad Z := \exp(\gamma I_1/\varepsilon) \quad \text{with } I_1 := I(x_1) < 0 \; . \tag{8.13}$$

This gives the rate for *each* $\gamma \geq 0$. Clearly $r_e = 0$ for $\gamma = 0$ (no diffusion), and $Z$ only contributes when $\gamma = O(\varepsilon)$. The result for lowest damping $r_e \approx \gamma|I_1|/\varepsilon$ is the same as in [6] ($\partial R/\partial\gamma = O(\varepsilon^0)$, see (8.14) below).

For the evaluation of $R$ it is sufficient to use the linear approximation $v_s(x) \approx \lambda_- x$

(so $I(x) \approx \lambda_- x^2/2$ and $x_1 \approx -\infty$), because $-g'(x) \to 2\delta(x)$ as $\varepsilon \to 0$. The result is just the well-known Kramers (not only Smoluchowski) rate factor

$$R = \omega/|\lambda_-|, \qquad (8.14)$$

as is easily inferred by use of the eigenvalue equation $\lambda^2 + \gamma\lambda - \omega^2 = 0$ for $\lambda_-$. Actually, (8.14) gives the double transition rate of Kramers. This shows that, given an arrival on the separatrix, both a transition and a return are equally probable (for smooth thresholds, and except for $\gamma = O(\varepsilon)$).

*Important:* No information about the inaccessible domain $x > 0$ has been used. The flaw of the existing approach is thus removed.

**Appendix :** *Asymptotic solutions of homogeneous backward equations.*

Consider a backward equation $L^+ q = 0$ with a diffusion term $\propto \varepsilon$, and a trial function $q^*(\vec{x}, \varepsilon)$ satisfying the (inhomogeneous) boundary condition for each $\varepsilon$. Clearly $-L^+ q^* := \theta(\vec{x}, \varepsilon)$ may be viewed as a source term for the errors, and $q^*$ is thus an asymptotic solution, if

$$\iint |\theta| dx\,dy \to 0 \text{ as } \varepsilon \to 0.$$

This admits persisting nonzero peaks of $\theta$ (with a width $\to 0$), as can be seen by the equation $(1+x)f_x + \varepsilon f_{xx} = 0$ ($x \geq 0$) for the probability integral $\Phi$: the well-known asymptotic solution $q^* = 1 - \exp(-x/\varepsilon)$ yields $\theta = (-x/\varepsilon)\exp(-x/\varepsilon)$ with $\max|\theta| = e^{-1}$ for each $\varepsilon > 0$, but with a width $\propto \varepsilon$.

*Application to the Kramers problem*

Insertion of (8.8) into (8.5) yields

$$\theta = -2g\varepsilon^{-1}\exp(2v_s\varepsilon^{-1}\eta)(\gamma v_s\eta + v_s'\eta^2) + (v_s+\eta)g'[1-\exp(2v_s\varepsilon^{-1}\eta)] \qquad (A.1)$$

($\eta := v - v_s \leq 0$, $x_1 \leq x \leq 0$). Clearly, for $\varepsilon > 0$

$$\theta = 0 \quad \text{on } \partial\Omega \ (\eta = 0), \text{ and for } x = 0 \ (\text{where } g = 0 = v_s). \qquad (A.2)$$

For $v_s\eta < 0$ the exponential function (times $\varepsilon^{-1}$) vanishes everywhere when $\varepsilon \to 0$, and thus the whole first term (recall that $0 \leq g \leq 1$). Moreover, by (8.7), $g'$ vanishes in the same sense wherever $U(x) < 0$, while at $x = 0$ is becomes $\delta(x)$, but the square bracket vanishes there. It is thus sufficient to show that $\theta$ remains bounded when $\varepsilon \to 0$. This is easily verified; we only note that the peak of the most critical term $-g\gamma(2v_s\varepsilon^{-1}\eta)\exp(2v_s\varepsilon^{-1}\eta)$ moves to $\partial\Omega$ when $\varepsilon \to 0$, without changing its height $g\gamma/e$, but with the width $\propto \varepsilon[2v_s(x)]^{-1}$; mind that $g(x)/v_s(x)$ exists at $x = 0$ and is $\propto \varepsilon^{-1/2}$ there. The integral of $|\theta|$ vanishes thus as $\varepsilon^{1/2}$.

**Concluding remarks**

1) For the exit problem it is required that the quasipotential $\phi$ exists and is differentiable at the preferred exit point. Whether possible singularities at different places can be relevant remains to be investigated.

2) A possible correction of $\phi$ by $\phi + \varepsilon\varphi$ multiplies the exit rate by $\exp(-\varphi_0)$, where $\varphi_0$ denotes $\varphi$ at the exit point (while $\varphi = 0$ at A). It is not difficult to show that $\varphi$ can be computed along the same characteristics as $\phi$.